\documentclass[12pt]{article}
\usepackage[centertags]{amsmath}
\usepackage{amsfonts}
\usepackage{amssymb}
\usepackage{amsthm}
\usepackage{amscd}% diagrams
\usepackage{color}
\theoremstyle{plain}
\newtheorem{thm}{Theorem}[section]

\newtheorem{lem}[thm]{Lemma}

\theoremstyle{definition}

\newtheorem{rem}[thm]{Remark}
\newtheorem{ex}[thm]{Example}
\newcommand{\om}{\omega}
\newcommand{\Om}{\Omega}
\renewcommand{\a}{\alpha}
\renewcommand{\b}{\beta}
\newcommand{\fee}{\varphi}

\newcommand{\2}{\rightarrow}
\newcommand{\duzhky}[1]{{\left ( #1 \right )}}
\newcommand{\bundlei}{{\mathcal I}}

\newcommand{\mrm}[1]{{\mathrm #1 } }
\newcommand{\Rmstar}{{ \duzhky{R_m}_*\, }}
\newcommand{\mt}[1]{ {\text{\tiny #1 }} }

\newcommand{\abs}[2]{ {|\, #1 |^{\, #2}} }
\newcommand{\ra}{\rightarrow}
\newcommand{\ssst}{\scriptscriptstyle}

\newcommand{\im}{{\text{\rmfamily\upshape Im} \,}}
\renewcommand{\Im}{{\im}}
\newcommand{\imH}{{\im\mathbb H}}
\renewcommand{\Re}{{\text{\rmfamily\upshape Re} \,}}
\newcommand{\Hstar}{{\mathbb H^*}}
\newcommand{\muo}{ {\mu ^{-1}(0)} }
\newcommand{\Tdot}{{T_{\boldsymbol .}\, }}
\DeclareMathOperator{\lspan}{\mathtt{span}}

\newcommand{\hK}{hyperK\"{a}hler }
\newcommand{\qK}{quaternionic K\"{a}hler }
\newcommand{\wrt}{with respect to }

\newcommand{\Kvf}{{Killing vector field }}
\newcommand{\kahler}{{K{\" a}hler }}
\begin{document}
%%% ====== BEGIN ======= BEGIN ======== BEGIN =====================

\title{HyperK{\"a}hler and quaternionic K{\"a}hler manifolds with $S^1$--\,symmetries}

\author{Andriy Haydys\thanks{e--mail: haydys@math.uni-bielefeld.de}\\ \textit{University of Bielefeld}}
\date{September 13, 2007}
%%\date{September 9, 2005}
\maketitle

\begin{abstract}
We study relations between quaternionic Riemannian manifolds admitting different types of
symmetries. We show that any hyperK\"ah\-ler manifold admitting \hK potential and
triholomorphic action of $S^1$ can be constructed from another \hK manifold (of lower
dimention) with an action of $S^1$ which fixes one complex structure and rotates the
other two and vice versa. We also study corresponding \qK manifolds equipped with a \qK
action of the circle. In particular we show that any positive \qK manifolds with
$S^1$--symmetry admits a \kahler metric on an open everywhere dense subset.
\end{abstract}

\section{Introduction}
In a short Berger's list \cite{Berger:55} of groups that can occur as holonomy group
of locally irreducible Riemannian manifold are in particular $Sp(n)$ and
$Sp(n)Sp(1)$, the first being a holonomy group of \hK manifold while the second of
\qK one. In other words, a Riemannian manifold $(M,g)$ is \hK if it admits three
covariantly constant complex structures $I_r,\ r=1,2,3$ with quaternionic
relations%
\begin{equation*}
I_r^2=-id,\qquad I_1I_2=-I_2I_1=I_3,
\end{equation*}
compatible with Riemannian structure: $g(I_r\cdot ,I_r\cdot)=g(\cdot ,\cdot)$.
Hitchin \cite{Hitchin:87} proved the following criteria of integrability of complex
structures: $I_r$ are covariantly constant if and only if 2--forms $\om_r(\cdot
,\cdot )=g(\cdot ,I_r\cdot)$ are symplectic (i.e. closed). It is convenient to
consider all three symplectic forms as a single 2-form with values in imaginary quaternions:%
\begin{equation*}
\om = \om_1 i+\om_2 j+\om_3 k.
\end{equation*}

In contrast to \hK manifolds, a \qK manifold $N$ admits complex structures (and
correspondingly 2--forms $\om_r$) only locally. Nevertheless the 4--form
$\Om=w_1\wedge\om_1 + w_2\wedge\om_2 +w_3\wedge\om_3$, called the fundamental
4--form, exists globally and determines the \qK structure. In this case the
integrability means that the fundamental 4--form is covariantly constant and this is
equivalent to $d\Om =0$ provided $\dim N\ge 12$. In dimension 4 \qK by definition
means  Einstein and self--dual.

%The first nontrivial example of \hK metric was one constructed by Eguchi and Hanson
%\cite{EguchiHanson:78} on $ T^*\mathbb{P}^1$ which was lately generalised by Calabi
%\cite{Calabi:79} to all $T^*\mathbb{P}^n$. In dimension four Gibbons and Hawking
%\cite{GibbonsHawking:78} described all \hK manifolds with $S^1$--symmetry (see example
%\ref{Ex_GibbonsHawking} for details).
%
%Another useful source of \hK manifolds is (finite- and infinite--dimensional)
%\hK reduction: if a group $G$ acts isometrically and triholomorphically on a
%\hK manifold $M$ with momentum map $\mu:M\2 \mathfrak g^*\otimes \imH$ and
%$a\in \mathfrak z\otimes\imH$ is its regular value, then $\mu^{-1}(a)/G$
%inherits \hK structure.

An important link between \hK and \qK geometries provides Swann
construction~\cite{Swann:91}. Suppose that the group $SU(2)\cong Sp(1)=$ $\left \{q\in
\mathbb H \mid |\, q|=1\right \}$ acts on $M$ isometrically and
\textit{permutes complex structures:}%
\begin{equation}\label{PermutingAction}
\mathcal L_{\zeta}\om=\left [ \zeta , \om \right ]\quad \Leftrightarrow \quad
\left (L_q\right )^*\om=q\om\bar q,
\end{equation}
where $\zeta\in \imH\cong \mathfrak {sp}(1)$. We also say, that the action of
$Sp(1)$ with the above property is \textit{permuting}. Swann shows that such
action can be extended to homothetic action of the whole $\Hstar=\mathbb
R_+^*\times Sp(1)$ if the vector field $IY_{I}$ is independent of a complex
structure $I$,
where $Y_I$ is a \Kvf of $S^1\subset Sp(1)$ which preserves $I$. In particular%
\begin{equation}\label{Eq_ConditionHstarAction}
I_1Y_1=I_2Y_2=I_3Y_3=-Y_0,
\end{equation}
where we put $Y_r=Y_{I_r}$ for short and a vector field $Y_0$ generates
\textit{homothetic} action of $\mathbb R_+^*\subset \mathbb H^*:\
\duzhky{L_r}^* g=r^2 g$. We will also call such $\Hstar$--action
\textit{permuting}. Under these circumstances $N=M/\Hstar$ has positive scalar
curvature and carries a \qK structure. On the other hand, for any \qK manifold
$N$ with positive scalar curvature Swann constructs  a \hK manifold $\mathcal
U(N)$ which enjoys permuting action of $\Hstar$. Such manifolds are also
distinguished by the property of carrying a \hK potential, i.e. function
$\rho:M\2\mathbb R$ which is K\"ahler potential simultaneously for each complex
structure.

%; this function is momentum map of any $S^1\subset Sp(1)$ \wrt the symplectic
%form which is preserved.

In this paper we study influence of $S^1$--symmetry on relations between different types
of quaternionic Riemannian geometries. We consider \hK manifolds $M$ with \hK potential
and additional triholomorphic and isometric action of $S^1$ and show that such manifolds
can be reconstructed from their \hK reductions $\tilde M$ with respect to a
\textit{nonzero} value of momentum map. The main result of this paper is
Theorem~\ref{Thm_HKmanifold} which describes $M$ as a certain bundle over $\tilde M$ with
the fibre $\Hstar$. Moreover, \hK structure is described quite explicitly which allows to
obtain not only existence results but metric and symplectic forms themselves. Dividing
such manifolds by $\Hstar$ as described in \cite{Swann:91} we obtain quaternionic
K\"ahler manifolds with $S^1$--symmetry (Theorem \ref{Thm_QKmflds}). In
Section~\ref{Section_IndeterminacyOfV} we describe  new examples of hyper-- and
quaternionic--Kahler manifolds making use of a certain freedom in choice of parameters of
the construction.

We also show in Theorem~\ref{Thm_KahlerStructure} that positive \qK manifold $N$
with $S^1$--symmetry admits a K\"ahler structure on an open everywhere dense
submanifold. The complex structure is a section of the structure bundle of $N$,
however the K\"ahler metric is different from the \qK one.

\section{$S^1$--Symmetry}

Let $M$ be a \hK manifold with $\Hstar$--action permuting complex structures.
Suppose also that $M$ admits a \hK action of $S^1$ (which we prefer to denote
$S_0^1$\index{Group!$S_0^1$} in order to distinguish this group from another one,
which will appear
later and is also isomorphic to $S^1$ ) with momentum map\index{Momentum map} $\mu: M\2 \imH$,%
\begin{equation*}
d\mu = -\imath_{K_0}\om,
\end{equation*}
where $K_0$ is a \Kvf of $S_0^1$. We assume also that these two actions
commute and that $\mu$ is $\Hstar$--equivariant:%
\begin{equation}\label{Mu_and_Hstar}
\mu\circ L_x=x\mu\bar x, \quad x\in \mathbb H^*.
\end{equation}
Fix an imaginary quaternion, say $i$, and consider corresponding level set
$P=\mu^{-1}(i).$  Since $xm=L_xm \in \mu^{-1}(xi\bar x)$ and $\mathbb H^*$ acts
transitively on $\imH\setminus\{0\}$ the map%
\begin{equation}\label{Mapf}
f:\mathbb H^*\times P\2 M\setminus\muo,\quad (x,m)\mapsto xm,
\end{equation}
is surjective. Notice that%
\begin{equation*}
M_0=M\setminus \mu^{-1}(0)
\end{equation*}
is open and everywhere dense submanifold of $M$.

However the map $f$ is not injective. Indeed $\mu^{-1}(i)=P$ inherits action of
$\mathrm{Stab}_i=S^1\subset \mathbb H^*$ (to which we now give a label
$S_r^1$\index{Group!$S_r^1$} ) and it follows that points $(x,m)$ and $(xz,\bar
zm),\ z\in S_r^1$ are mapped into the same point $xm$ in $M_0$. Thus the manifold
$M_0$ can be described as $\Hstar\times_{S_r^1}P.$ Now the challenge is to express
the \hK structure of $M_0$ in terms of its "components" $\Hstar$ and $P$.

While the first "component" $\Hstar$ is very simple the second one needs to be
understood more deeply for the future purposes.

\subsection{Induced structure on P}\label{Subsection_IndStrOnP}

It follows from $\Hstar$--equivariancy of $\mu$ that each nonzero imaginary quaternion is
a regular value of $\mu$. Assuming that $S^1_0$ acts freely on $P$, we see that $\tilde
M:=P/S_0^1$ is just a \hK reduction of $M$ and therefore is itself a \hK manifold. Thus
$P$ can be thought of as $S_0^1$--principal bundle over $\tilde M$. Moreover it comes
equipped with a connection\footnote{Notice that we identify the Lie algebra of $S^1$
with $\mathbb R$, not with $i\mathbb R$}, namely%
\begin{equation*}
\xi(\cdot)=v g(K_0,\cdot)\in \Omega^1(P),
\end{equation*}
where $v^{-1}=g(K_0,K_0),\ v:\tilde M\2 \mathbb R_{>0}$. Notice that the induced
metric $\tilde g$ on $\tilde M$, the connection $\xi$ and the function $v$ together
determine the metric on $P$ since $T_{\boldsymbol .}\, P\cong \mathbb R K_0\oplus
T_{\boldsymbol .}\, \tilde M:$
\begin{equation}\label{MetricOnP}
g_{\mt P}=\tilde g + v^{-1}\xi^2.
\end{equation}
The connection $\xi$ defines a horizontal lift $\hat {\mrm u}\in TP$ of a tangent
vector $\mrm u\in T\tilde M$.

As we have already remarked $P$ inherits the action of $S_r^1$, which descends to
$\tilde M$. The latter action has a nice property (inherited from $M$) of fixing
complex structure $I_1$ and rotating the plane spanned by $I_2$ and $I_3$. Denote by
$K_r$ a Killing vector field of $S_r^1$--action on $\tilde M$
and by $w$ the squared norm of $K_r$:%
\begin{equation*}
w:\tilde M\ra \mathbb R_{>0},\quad w=\|K_r\|^2.
\end{equation*}
Below we will also use a quaternion--valued 1--form $\eta$ generated by $K_r$:%
\begin{equation}\label{Eq_FormEta}
\eta=\imath_{K_r}\tilde g + \imath_{K_r}\tilde\om \in\Om^1(\tilde M; \mathbb H).
\end{equation}
Further, recall that $Y_1$ is the \Kvf of the $S_r^1$--action on $P$. Then $Y_1$ and
$K_r$ are related as follows. First observe that $\Tdot M = \Tdot P\oplus\mathbb
RI_1K_0 \oplus\mathbb
RI_2K_0\oplus\mathbb RI_3K_0$ and one also has%
\begin{equation*}
\Tdot P=\mathrm{Ker}\, \mu_*,\quad \mu_* I_1K_0=v^{-1}i,\ \mu_* I_2K_0=v^{-1}j\
\text{ and}\quad \mu_* I_3K_0=v^{-1}k.
\end{equation*}
Now taking $x=\exp (it)$ in formula (\ref{Mu_and_Hstar}) and differentiating \wrt
$t$ one obtains that the formula $Y_1=\hat K_r+aK_0$ holds on $P$. The same argument
gives that $\mu_*Y_0=2i$ or in other words $Y_0=\hat Y'+ bK_0+2vI_1K_0$. It follows
from the equation $I_1Y_0=Y_1$ that $Y'=-I_1K_r,\ b=0,\ a=-2v.$ Summing up we
obtain%
\begin{equation}\label{Y_andK_r}
\begin{gathered}
Y_0=-I_1\hat K_r - 2vI_1K_0,\quad%
Y_1=\hat K_r + 2vK_0,\qquad\ \ \\%
Y_2=I_3\hat K_r + 2vI_3K_0,\quad%
Y_3=-I_2\hat K_r - 2vI_2K_0.%%
\end{gathered}
\end{equation}

\begin{rem}\label{Rem_EquivariancyOfXi}
Since actions of $S_0^1$ and $S_r^1\subset \Hstar$ commute, it follows that the
connection $\xi$ enjoys additional property of being $S_r^1$ invariant. On
infinitesimal level this means that $0=\mathcal
L_{Y_1}\xi=\imath_{Y_1}d\xi+d\imath_{Y_1}\xi=\imath_{K_r}F_\xi+2dv,$ where
$F_\xi\in\Om^1(\tilde M)$ denotes the curvature form of $\xi$. Thus, invariance of
$\xi$ \wrt action of $S_r^1$
on $P$ is equivalent to%
\begin{equation}\label{Eq_invariancy_of_xi}
\imath_{K_r}F_\xi+2dv=0.
\end{equation}
Note also that the function $v$ is $S_r^1$--invariant by the same reason.
\end{rem}

\subsection{Metric}\label{Section_Metric}

Since $M$ is a Riemannian manifold the map $f$, defined by (\ref{Mapf}), induces a
metric $f^*g$ on $\Hstar\times P$. Notice that since $f$ is not injective, this
metric degenerates on tangent vectors to fibres. Our next aim is to calculate $f^*g$
explicitly in terms of tensors on $\Hstar$ and $\tilde M$ as well as the connection
$\xi$ and function $v.$

Let $(x,m)\in \Hstar\times P$ and $(h_1  ,\mathrm v_1),\ (h_2 , \mathrm v_2)\in
T_x\Hstar\times T_m P$. Put also $\a = x^{-1}h_1, \b = x^{-1}h_2\in T_1\Hstar$ and
denote by $\mathsf Y_\a$ and $\mathsf Y_\b$ the Killing vector fields of
$\Hstar$--action at the point $m$ corresponding to the Lie algebra elements $\a$ and
$\b$. Obviously $(\mathsf Y_1, \mathsf Y_i, \mathsf Y_j, \mathsf Y_k)=(Y_0, Y_1,
Y_2, Y_3)$.
Further, one has%
\begin{align*}
f^*g \bigl ( (h_1, \mathrm v_1), (h_2, \mathrm v_2) \bigr )&=%
g\bigl ( \duzhky{L_x}_*\duzhky{\mathsf Y_\a +\mathrm v_1},\,%
    \duzhky{L_x}_*\duzhky{\mathsf Y_\b +\mathrm v_2} \bigr )\\%
&=|\, x|^{\, 2}g\bigl ( \mathsf Y_\a +\mathrm v_1 ,\, %
    \mathsf Y_\b +\mathrm v_2 \bigr ).
\end{align*}
Thus we see that essentially the following three terms have to be computed:
$g\duzhky{ \mathsf Y_\a , \mathsf Y_\b },\ g\duzhky{ \mathsf Y_\a , \mathrm v }$ and
$g\duzhky{\mathrm v_1,\mathrm v_2}$.

\textbf{The first term.} Since relation (\ref{Eq_ConditionHstarAction}) holds, we
get
\begin{equation*}
g \duzhky{ \mathsf Y_\a , \mathsf Y_\b }=%
  g\duzhky{\sum_{r=0}^3 \a_rY_r\ ,\sum_{r=0}^3\b_rY_r }=%
 g\duzhky{Y_0,Y_0}\Re\duzhky{\a\bar\b}.
\end{equation*}
Recall that $w$ denotes the squared norm of $K_r$ and therefore it follows from
(\ref{Y_andK_r}) that $g(Y_0,Y_0)=w+4v^2v^{-1}=4v+w.$ So finally we have%
\begin{equation*}
g\duzhky{ \mathsf Y_\a , \mathsf Y_\b }= \duzhky{4v+w}\Re\duzhky{\a\bar\b}.
\end{equation*}

\textbf{The second term.} First decompose $\mathrm v$ into horizontal and vertical
parts: $\mathrm v= \hat{\mathrm v}'+ \xi (\mathrm v)K_0$. Taking into
account formulae (\ref{Y_andK_r}) again, one obtains %
\begin{equation*}
\begin{split}
g (&\mathsf Y_\b , \mathrm v )=%
g\duzhky{\mathsf Y_\b , \mathrm v'}+ 2\b_1\xi (\mathrm v)\\%
&\quad =\tilde g\duzhky{-\b_0I_1K_r+\b_1K_r+\b_2I_3K_r-\b_3I_2K_r,
\mathrm v'} + 2\b_1\xi (\mathrm v)\\%
&\quad =\b_0\tilde\om_1(K_r, \mathrm v')+ \b_1\tilde g(K_r, \mathrm v')-
\b_2\tilde\om_3(K_r, \mathrm v')+ \b_3\tilde\om_2(K_r, \mathrm v') + 2\b_1\xi
(\mathrm v).
\end{split}
\end{equation*}
Slightly abusing notations, we also use the letter $\eta$ for the pull--back of the
form (\ref{Eq_FormEta}) to $P$. Then the above formula can be written in a more
compact form:
\begin{equation*}
g\duzhky{\duzhky{R_m}_*\b , \mathrm v} =%
 -\Re\bigl (2\b i\xi(\mathrm v)+\b i\eta(\mathrm v)\bigr ).
\end{equation*}

\textbf{The third term.} This has been already computed and is given by
(\ref{MetricOnP}).

\begin{rem} Below we follow conventions of \cite{GreubAO:72}. In
particular, if $\zeta_1$ and $\zeta_2$ are (quaternion--valued)
1--forms, then%
\begin{equation*}
\begin{gathered}
\duzhky{\zeta_1\odot\zeta_2}\duzhky{v_1,v_2}=\zeta_1(v_1)\zeta_2(v_2)+\zeta_1(v_2)\zeta_2(v_1),\\
\duzhky{\zeta_1\wedge\zeta_2}\duzhky{v_1,v_2}=\zeta_1(v_1)\zeta_2(v_2)-\zeta_1(v_2)\zeta_2(v_1).
\end{gathered}
\end{equation*}
\end{rem}

Now, recalling that $\a$ and  $\b$ contain shift by %
$x^{-1}=\bar x/|\, x|^{\, 2}$ we obtain a final form of the metric:%
\begin{equation}\label{Metric}
f^*g= (4v+w)\Re\, dx\otimes d\bar x -\Re\duzhky{\bar xdxi\odot (2\xi +\eta ) } +|\,
x|^{\, 2}\duzhky{\tilde g + v^{-1}\xi^2}.
\end{equation}

%---END---\subsection{Metric.}

\subsection{Symplectic forms}

In this section we will describe symplectic forms in the similar manner as we did
with the metric above.

The pull-back of $\om$ can be written as%
\begin{align*}
f^*\om ((h_1, \mrm v_1),(h_2,\mrm v_2) )&=\om\duzhky{ \duzhky{L_x}_*\duzhky{\mathsf Y_\a +\mrm v_1},%
    \duzhky{L_x}_*\duzhky{\mathsf Y_\b +\mrm v_2}  }\\%
&=x\om\duzhky{ \mathsf Y_\a +\mrm v_1,\mathsf Y_\b +\mrm v_2 } \bar x,
\end{align*}
where $\a$ and $\b$ are the same as in Section \ref{Section_Metric}. Therefore we
have to compute three terms analogous to those, which appear in the metric
computation.

\textbf{The first term.} The computation is similar to the one above:%
\begin{align*}
\om (&\mathsf Y_\a,\mathsf Y_\b)=\\%
    &= ig\duzhky{\a_0Y_0+\a_1Y_1+\a_2Y_2+\a_3Y_3 ,
        \, \b_0Y_1-\b_1Y_0+\b_2Y_3-\b_3Y_2 }\\%
&+jg\duzhky{\a_0Y_0+\a_1Y_1+\a_2Y_2+\a_3Y_3 ,
        \, \b_0Y_2-\b_1Y_3-\b_2Y_0+\b_3Y_1 }\\%
&+kg\duzhky{\a_0Y_0+\a_1Y_1+\a_2Y_2+\a_3Y_3 ,
        \, \b_0Y_3+\b_1Y_2-\b_2Y_1-\b_3Y_0 }\\%
&=g(Y_0,Y_0)\left (\, i\duzhky{-\a_0\b_1 +\a_1\b_0-\a_2\b_3+\a_3\b_2}\right.\\%
     &+ \left. j\duzhky{-\a_0\b_2+\a_2\b_0 + \a_1\b_3-\a_3\b_1}+%
   k\duzhky{-\a_0\b_3+\a_3\b_0 - \a_1\b_2+\a_2\b_1} \right )\\%
&=(4v+w)\Im\duzhky{\a\bar\b}.
\end{align*}

\textbf{The second term.} Decomposing $\mrm v$ into horizontal
$\hat{\mrm v}'$ and vertical $\xi(\mrm v) K_0$ parts one obtains:%
\begin{align*}
&\om \duzhky{ \mathsf Y_\a , \mrm v }=%
    i\duzhky{-2\a_0\xi(\mrm v)+\om_1\duzhky{\Rmstar\a ,\hat{\mrm v}'} }\\%
    &+j\duzhky{ -2\a_3\xi(\mrm v)+\om_2\duzhky{\Rmstar\a ,\hat{\mrm v}'} }%
+k\duzhky{ -2\a_2\xi(\mrm v)+\om_3\duzhky{\Rmstar\a ,\hat{\mrm v}'} }\\%
&=i\duzhky{-2\a_0\xi(\mrm v)- \a_0\tilde g(K_r, \mrm v') +\a_1\tilde\om_1(K_r,%
    \mrm v') - \a_2\tilde\om_2 (K_r, \mrm v') - \a_3\tilde\om_3 (K_r, \mrm v') }\\%
&+j\duzhky{-2\a_3\xi(\mrm v)+ \a_0\tilde\om_3(K_r, \mrm v') +\a_1\tilde\om_2(K_r,%
    \mrm v') + \a_2\tilde\om_1 (K_r, \mrm v') - \a_3\tilde g(K_r, \mrm v') }\\%
&+k\duzhky{-2\a_2\xi(\mrm v)- \a_0\tilde\om_2(K_r, \mrm v') +\a_1\tilde\om_3(K_r,%
    \mrm v') + \a_2\tilde g(K_r, \mrm v') + \a_3\tilde\om_1(K_r, \mrm v') }\\%
&=-2\im\duzhky{\a i}\xi(\mrm v) - \im\duzhky{\a i\eta(\mrm v)} .
\end{align*}

\textbf{The third term.} It is easy to see that %
\begin{equation*}
\om (\mrm v_1 , \mrm v_2)=\om (\mrm v_1' , \mrm v_2')=%
    \tilde\om (\mrm v_1 , \mrm v_2),
\end{equation*}
where the pull--back is also implied.

\medskip

Thus, recalling that $\a = x^{-1}h_1=|x|^{-2}\bar x h_1$,
the $\imH$--valued form $\fee = f^*\om$ can be written as%
\begin{equation}\label{Symplectic_form}
\fee = \frac{4v+w}2 dx\wedge d\bar x + x\tilde\om\bar x -%
2\im\duzhky{dxi\bar x}\wedge\xi - \im\duzhky{dxi\wedge\eta\bar x}.
\end{equation}

%---END---\subsection{Symplectic forms.}

\subsection{Inverse Construction}\label{Section_ReverseConstruction}

Now we can look on the above considerations in reverse order in the following sense.
Suppose $\tilde M$ is a \hK manifold with metric $\tilde g$ and \hK structure
$\tilde\om$. Further, a group $S_r^1$ acts on $\tilde M$ preserving complex
structure $I_1$ and rotating $I_2$ and $I_3$ in the sense
$\duzhky{L_z}^*\tilde\om=z\tilde\om\bar z,\ z\in S_r^1$. Pick an $S_0^1$--principal
bundle $P$ with a connection $\xi$ and extend the action of $S_r^1$ to $P$ such that
it commutes with $S_0^1$ (at least locally such extension always exists).

Consider further a manifold $M_0=\Hstar\times_{S_r^1}P$. We would like to define a
metric $g$ and \hK structure $\om$ on $M_0$ such that their pull--backs to
$\Hstar\times P$ are given by formulae (\ref{Metric}) and (\ref{Symplectic_form})
respectively. The first thing to show is that these expressions define invariant and
basic tensors on $\Hstar\times P$. One can easily check that both tensors are
invariant provided $\xi$ is $S_r^1$--invariant (see also remark
\ref{Rem_InvariancyOfXi}). Let $\chi$ be a \Kvf of the $S_r^1$--action on
$\Hstar\times P$. It follows that $\chi=K^* -Y_1$, where $K^*$ is a \Kvf of the
$S_r^1$--action on $\Hstar$ by right multiplication, i.e. $dx(K^*)=xi$. Then the
equalities $\imath_\chi g=0,\ \imath_\chi\fee =0$ can be checked directly. For
example, the last one
follows from the following computation:%
\begin{align*}
\duzhky{\imath_\chi\fee}\duzhky{\a,\mrm v}=&%
    \frac 12 (4v+w)\duzhky{xi\bar\a -\a\overline{xi} } - %
    x\tilde\om\duzhky{K_r,\mrm v}\bar x\\%
&-2\im\duzhky{xii\bar x\xi \duzhky{\mrm v}+\a i\bar x\, 2v} -%
    \im\duzhky{ xii\eta\duzhky{\mrm v}\bar x +\a i\eta\duzhky{K_r}\bar x }\\%
&=(4v+w)\im\duzhky{xi\bar\a}-x\tilde\om\duzhky{K_r,\mrm v}\bar x -2\cdot 0 \\%
    &-4v\im\duzhky{xi\bar\a} +x\tilde\om\duzhky{K_r,\mrm v}\bar x - w\im{xi\bar\a} \\%
&=0.
\end{align*}

The next question is wether the 2--form $\om\in\Om^2(M_0; \Im\mathbb H)$ is closed.
As we have seen, the pull--back $\fee$ of $\om$ to $\Hstar\times P$ is basic and
therefore this is equivalent to $\fee$ being closed. Now $d\fee$ is a
quaternion--valued 3--form on $\Hstar\times P$ and by the K\"unneth formula
$\Om^3(\Hstar\times P; \Im\mathbb H)\cong \bigoplus_{l=0}^3\Om^l(\Hstar; \Im \mathbb
H)\otimes\Om^{3-l}(P; \Im \mathbb H)$. Thus $d\fee$ decomposes in 4
components: $d\fee = \sum_{l=0}^3 (d\fee )_{(l,\, 3-l)}$, %
$(d\fee)_{(l,\, 3-l)}\in \Om^l(\Hstar; \Im \mathbb H)\otimes\Om^{3-l}(P; \Im \mathbb
H)$.  It is easy to see that $(d\fee)_{(0,\, 3)}$ and $(d\fee)_{(3,\, 0)}$ vanish
identically and it remains to compute the remaining two components of $d\fee$.

It follows directly from the expression for $\fee$ that%
\begin{align*}
(d\fee)_{(1,\, 2)}&= dx\wedge\tilde\om\bar x+  x\tilde\om\wedge d\bar x +%
    2\Im\duzhky{dxi\bar x}\wedge F_\xi +\Im\duzhky{dxi\wedge d\eta\bar x}\\%
&=\Im\duzhky{dx\wedge\duzhky{2\tilde\om +2iF_\xi +id\eta}\bar x}
\end{align*}
and this vanishes iff%
\begin{equation*}
-2i\tilde\om + 2F_\xi+d\eta = 0.
\end{equation*}

By the Cartan formula $[i,\tilde\om]=\mathcal L_{K_r}\tilde \om =
d\duzhky{\imath_{K_r}\tilde\om}$. But then the above equation can be rewritten as
$2F_\xi=-d\duzhky{\imath_{K_r}\tilde g} -d\duzhky{\imath_{K_r}\tilde\om}
+2i\tilde\om=-d\duzhky{\imath_{K_r}\tilde g} -2\tilde\om_1$. Thus the vanishing of
$(d\fee)_{(1,\, 2)}$ is equivalent to%
\begin{equation}\label{Eq_F}
F_\xi=-\frac 12\, d\duzhky{\imath_{K_r}\tilde g} -\tilde\om_1.
\end{equation}
For the other nontrivial component of $d\fee$ one obtains%
\begin{align*}
(d\fee)_{(2,\, 1)}&=\frac 12 (4dv+dw)\wedge dx\wedge d\bar x +2\Im\duzhky{dxi\wedge d\bar x}\wedge\xi-%
    \Im\duzhky{dxi\wedge\eta\wedge d\bar x}\\%
&=\frac 12 \Im\duzhky{dx\wedge\duzhky{-\duzhky{4dv+dw}-4i\xi-2i\eta }\wedge d\bar
x}.
\end{align*}

Suppose $\theta$ is a quaternion--valued 1--form on $\tilde M$ and consider the
equation $\ \Im\duzhky{dx\wedge\theta\wedge d\bar x}=0$ on $\Hstar\times \tilde M$,
which turns out to be equivalent to $\Re\,\theta=0$. Indeed, $dx\wedge\theta\wedge
d\bar x=-\duzhky{\Re\, \theta}\wedge dx\wedge d\bar x + dx\wedge\Im\, \theta\wedge
d\bar x$ and the last summand is real--valued: $\overline{ dx\wedge\Im\,
\theta\wedge d\bar x }= (-1)\, dx\wedge\overline{\Im\theta}\wedge d\bar x =
dx\wedge\Im\, \theta\wedge d\bar x.$

Therefore $(d\fee)_{(2,\, 1)}$ vanishes iff%
\begin{equation}\label{Eq_dv}
4dv+dw=2\imath_{K_r}\tilde\om_1.
\end{equation}

Thus, the 2--form $\fee$ descends to a closed form on $M_0=\Hstar\times_{S_r^1}P$ if
and only if the three equations are satisfied : (\ref{Eq_invariancy_of_xi}),
(\ref{Eq_F}) and (\ref{Eq_dv}). But the last equation follows from the first two.
Indeed, since $S_r^1$ acts isometrically we have $0=\mathcal L_{K_r}\duzhky{
\imath_{K_r}\tilde g } = \imath_{K_r}\, d\duzhky{\imath_{K_r}\tilde g} +
d\duzhky{\imath_{K_r}\imath_{K_r}\tilde g}$ which means $\imath_{K_r}\,
d\duzhky{\imath_{K_r}\tilde g} = -dw$. Now taking the operator $\imath_{K_r}$ of
both sides of equation (\ref{Eq_F}) and using (\ref{Eq_invariancy_of_xi}) we obtain
equation (\ref{Eq_dv}).

It was first remarked in \cite{HitchinAO:87} that a \hK manifold with an
$S^1$--action which preserves one complex structure and permutes the other two has a
\textit{K\"ahler} potential. Since our conventions slightly differ we reproduce this
simple computation.

Let $\tilde\rho :\tilde M\2\mathbb R$ be a momentum map of $S_r^1$, i.e. a
solution of the equation%
\begin{equation*}
d\tilde\rho = -\imath_{K_r}\tilde\om_1.
\end{equation*}
On the one hand we have $d\duzhky{I_2^*d\tilde\rho}=d\, i\duzhky{\partial_2 -
\bar\partial_2}\tilde\rho= -2i\partial_2\bar\partial_2\tilde\rho$. But on the other
hand $I_2^*\imath_{K_r}\tilde\om_1 = \imath_{K_r}\tilde\om_3$ and therefore
$-d\duzhky{I_2^*d\tilde\rho}=d\duzhky{I_2^*\imath_{K_r}\tilde\om_1}$
$=d\duzhky{\imath_{K_r}\tilde\om_3}=\mathcal L_{K_r}\tilde\om_3 =2\tilde\om_2$.
Putting this together we obtain that $\tilde\rho$ satisfies%
\begin{equation*}
i\partial_2\bar\partial_2\tilde\rho = \tilde\om_2
\end{equation*}
or, in other words, $\tilde\rho$ is a \textit{K\"ahler} potential for $\tilde\om_2$.
It is clear that $\tilde\rho$ is also a K\"ahler potential for $\tilde\om_3$ since
these forms are not distinguished by the $S_r^1$--action. However $\tilde\rho$ needs
not to be a K\"ahler potential for $\tilde\om_1$.

Now if we remark that $I_1^*d\tilde\rho=\imath_{K_r}\tilde g$ and consequently%
\begin{equation*}
-2i\partial_1\bar\partial_1\tilde\rho = d\duzhky{\imath_{K_r}\tilde g},
\end{equation*}
then equation (\ref{Eq_F}) can be written in a particularly nice form:
$F_\xi=i\partial_1\bar\partial_1\tilde\rho - \tilde\om_1$, i.e. the function
$\tilde\rho$ is a \hK potential iff $F_\xi=0$.

Now we can find the function $v$ (up to a constant) from equation (\ref{Eq_dv})
(or, equivalently, from (\ref{Eq_invariancy_of_xi})):%
\begin{equation}\label{Function_v}
v =-\frac{w+2\tilde\rho}4.
\end{equation}

Therefore the following theorem is essentially proven.

\begin{thm}\label{Thm_HKmanifold}
Let a group $S_r^1$ act isometrically on a \hK manifold $\tilde M$ such that
$\duzhky{L_z}^*\tilde\om = z\tilde\om\bar z$. Further, let $P\2 \tilde M$ be an
$S_0^1$--principal bundle with a connection $\xi\in\Om^1(P)$. Suppose also, that the
function $v$ defined by formula (\ref{Function_v}) is everywhere positive, where $w$
denotes squared norm of the \Kvf $K_r$ of $S_r^1$, while $\tilde\rho$ is its
momentum map. Extend the action of $S_r^1$ to $P$ such that it commutes with the
action of $S_0^1$. Then (\ref{Metric}) and (\ref{Symplectic_form}) define a \hK
structure on $M_0=\mathcal H(\tilde M)=\Hstar\times_{S_r^1} P$ if and only if
\begin{equation}\label{Eq_Main}
F_\xi=i\,\partial_1\bar\partial_1\tilde\rho - \tilde\om_1\, .
\end{equation}
Furthermore the left action of $\Hstar$ induces a transitive action on the 2--sphere
of complex structures and therefore $\mathcal H(\tilde M)$ has a \hK
potential%
\begin{equation}\label{Eq_hKPotential}
\rho= -\frac {4v+w}2\abs {x}2.
\end{equation}
Finally, for any \hK manifold $M$ with permuting action of $\Hstar$ and
triholomorphic one of $S^1$, the open everywhere dense submanifold $M_0=M\setminus
\mu^{-1}(0)$ can be obtained as $\mathcal H(\tilde M)$, where $\tilde M$ is as
above.
\end{thm}

\begin{proof}
It remains to show that the symmetric tensor given by formula (\ref{Metric})
provides a non--negative bilinear form at any point of the tangent space to
$\Hstar\times P$ as well as to prove formula (\ref{Eq_hKPotential}) for the \hK
potential.

First we have a decomposition $T\duzhky{\Hstar\times P}=T\mathbb H\oplus TP
=T\mathbb H\,\oplus\, \mathbb RK_0\oplus \pi^*T\tilde M$. Further decompose $T\tilde
M$ as $\lspan( K_r, I_1K_r, I_2K_r, I_3K_r)\oplus E$, where $E$
denotes the orthogonal complement. Thus we have%
\begin{equation*}
T\duzhky{\Hstar\times P}= T\mathbb H \oplus \mathbb RK_0\oplus%
\pi^*\mathtt{span}\duzhky{K_r, I_1K_r, I_2K_r, I_3K_r}\oplus \pi^* E,
\end{equation*}
and we can write a tangent vector as $\mathrm w=\mrm w^* + aK_0 + \b K_r + \mathrm v$,
where $a$ is a real number and $\b$ is a quaternion\footnote{Any tangent space of a \hK
manifold carries an action of $\mathbb H$. In particular if $\b=\b_0+\b_1
i+\b_2j+\b_3k\in\mathbb H, \ $ we write $\b K_r$ instead of $\b_0K_r+\sum_{l=1}^3\b_l
I_lK_r$ for the sake of brevity}.
If $dx(\mrm w^*)=\a\in\mathbb H$, then%
\begin{align*}
g(\mathrm w, \mathrm w)%
&=(4v+w)\abs\a 2 +\abs x2\duzhky{\abs \b 2 + \|\mathrm v\|^2 + \frac {a^2}v} -%
    2\Re\duzhky{\bar x\a i\duzhky{2a+w\b}}\\%
&=4v\abs \a 2 - 4a\Re\duzhky{\bar x\a i} +\abs x2 v^{-1} a^2 \\%
&%\hskip4cm
+w\duzhky{\abs \a 2 - 2\Re\duzhky{\bar x\a i\b}+\abs{\b x}2}%
+\abs x2 \|\mathrm v\|^2 \\%
&=\left | 2\sqrt v\a i - \frac {ax}{\sqrt v} \right |^2 +%
    w\bigl | \a i-\b\bar x \bigr |^2 + \abs x2 \|\mathrm v\|^2\ge 0.
\end{align*}
Further, it was shown in \cite{Galicki:93} that if a \hK manifold $M$ admits a permuting
$\Hstar$--action, then the squared norm of any \Kvf generating this action is a \hK
potential (up to a constant\footnote{the minus sign appears because of different sign
convention in the definition of \hK potential} $-2$). Now the permuting action of
$\Hstar$ on $\mathcal H(\tilde M)$ is induced by the left multiplication on the first
component of $\Hstar\times P$. In particular, the \Kvf of $\mathbb R^*\subset\mathbb H^*$
is the vector field $\mrm w^*$ s.t. $dx(\mrm w^*)=x$. Its squared norm (multiplied by
$-1/2$) \wrt metric (\ref{Metric}) is exactly the right--hand side of
(\ref{Eq_hKPotential}).\qed
\end{proof}
\begin{rem}
It is easy to see that the \hK reduction of $\mathcal H (\tilde M)$ by $S_0^1$ is $\tilde
M$ (certainly not a surprise in view of Section~\ref{Subsection_IndStrOnP}). Thus the
construction $\mathcal H(\cdot )$ may be regarded as a kind of "\hK induction", i.e. an
inverse construction to the \hK reduction.
\end{rem}

\begin{rem}\label{Rem_InvariancyOfXi}
From equation (\ref{Eq_Main}) we have that the Chern class of $P$ is $-\frac
1{2\pi}[\,\tilde\om_1]$. It follows that $\frac 1{2\pi}[\,\tilde\om_1]$ must represent an
inegral cohomology class and in this case equation (\ref{Eq_Main}) does have a solution.
Furthermore, any solution $\xi$ is automatically $S_r^1$--invariant and the Killing
vector field $Y_1$ of $S_r^1$--action on $P$ satisfies $Y_1=\hat K_r + 2vK_0$. Indeed, as
we have seen the right hand side of equation (\ref{Eq_Main}) may be written in the form
$-\frac 12d\duzhky{\imath_{K_r}\tilde g} -\tilde\om_1$ and this immediately implies
$S_r^1$--invariancy of $\xi$. Further, we may decompose $Y_1$ on the horizontal and
vertical parts: $Y_1=\hat K_r + aK_0$. Then by Remark~\ref{Rem_EquivariancyOfXi} we have
$\imath_{K_r}F_\xi+da=0$. On the other hand equation (\ref{Eq_Main}) implies
$\imath_{K_r}F_\xi=-2dv$ and the statement follows.

It is worth pointing out that equality $a=2v$ holds only up to a constant. This
phenomenon will be discussed in details in Section~\ref{Section_IndeterminacyOfV}. At
this point we will ignore this subtlety implying that a constant is chosen properly, i.e.
such that equation $Y_1=\hat K_r + 2vK_0$ holds.
\end{rem}

%\begin{rem}[Complex structures]
%We would like to indicate how one can describe complex structures of $\mathcal
%H(\tilde M)$. It is convenient for a while to index complex structures of a \hK
%manifold by imaginary quaternions $q$ of unitary length. Let $I_q$ be one of the
%complex structures of $M$ and $(h,\mrm v)\in T\Hstar\times TP$.
%Consider the equation%
%\begin{equation}\label{Eq_ComplexStructure}
%f_* (h_1, \mrm v_1)=I_q f_*(h,\mrm v),
%\end{equation}
%where the map $f$ is as defined in (\ref{Mapf}).  Its solution $(h_1,\mrm
%v_1)=J_q(h,\mrm v)$ is defined only up to the \Kvf $\chi$; although the operator
%$J_q$ is not defined on $T(\Hstar\times P)$, it is well--defined on the factorspace
%$T(\Hstar\times P)/\lspan(\chi)$. Since we are free to add any vector of $\lspan
%(\chi)$ to both $(h,\mrm v)$ and $(h_1,\mrm v_1)$ it can be assumed that $\xi (\mrm
%v)=0=\xi (\mrm v_1)$. Under these conventions equation (\ref{Eq_ComplexStructure})
%has a unique solution $J_q(h,\mrm v)=(qh, I_{x^{-1}qx}\,\mrm v)$.
%\end{rem}

\begin{rem}\label{Rem_CaseOfHKPotential}
Suppose that the $S_r^1$--action is induced by a permuting $\, \Hstar$--action (and
standard inclusion $S_r^1\subset\Hstar$), or equivalently, the momentum map of the
$S_r^1$--action is not only K\"ahler potential but also {\hK}\cite{Swann:91}. It
follows from equation (\ref{Eq_Main}) that the bundle $P$ is flat and we can take it
to be trivial so that topologically $\mathcal H(\tilde M)=\Hstar\times \tilde M$.
Moreover, it follows from the proof of the theorem that $v$ is constant so that we
may put $v=1$. This determines a metric and symplectic forms.

Further, it turns out that in this case $\mathcal H(\tilde M)$ is \textit{isometric}
to $\Hstar\times \tilde M$ with its product metric. Indeed, direct computation shows
that the map $\mathcal H(\tilde M)\2\Hstar\times \tilde M,\ (x,m)\mapsto (x, xm)$ is
an isometry.
\end{rem}

%---END---\subsection{Reverse Construction.}

\subsection{Quaternionic Flip}\label{Section_QuaternionicFlip}

In the previous section for any \hK manifold $\tilde M$ with a certain
$S^1$--symmetry we have constructed another \hK manifold $M_0=\mathcal H(\tilde M)$
with \hK potential. Then Swann's results \cite{Swann:91} imply that the ma\-ni\-fold
$N_0=M_0/\Hstar=P/S_r^1$ is quaternionic K\"ahler. In this section we will describe
its \qK structure.

First notice that in order to obtain \qK structure on $N_0$ we have to consider a
\textit{riemannian} version of the quotient $M_0/\Hstar$, that is to pick a level
set of a \hK potential and divide it by the group $Sp(1)$; in this case we may view
complex structures of $N_0$ as induced by those of $M_0$ on $\lspan(Y_0, Y_1, Y_2,
Y_3)^\bot\subset TM_0$.

Let us again return to the viewpoint of Section \ref{Subsection_IndStrOnP}, i.e.
$P=\mu^{-1}(i)\subset M$ and let $\lambda=(4v+w)^{-1/2}$. Since the
restriction of the \hK potential $\rho$ to $P$ equals $-(4v+w)/2$, a map%
\begin{equation}\label{Eq_Isomorphism_l}
l:\, p\mapsto \lambda(p)\cdot p=L_{\lambda(p)}\, p,\quad p\in P,
\end{equation}
is a diffeomorphism between $P$ and $Q=\rho^{-1}(-1/2)\cap \mu_c^{-1}(0)\cap
\{\mu_1>0 \}$, where $\mu_c=\mu_2+i\mu_3$. Thus our next aim is to compute the
tensors $g(pr\circ l_*\,\cdot\, , pr\circ l_*\,\cdot)$ and $\om(pr\circ l_*\,\cdot\,
, pr\circ l_*\,\cdot)$, where $pr$ means a projection onto $\lspan(Y_0, Y_1, Y_2,
Y_3)^\bot$.

First we may decompose a vector $\mrm u\in T_pP\subset T_pM$ as $\mrm
u'+\sum_{l=0}^3a_lY_l$. The coefficients $a_l$ can be found from the following relations:%
\begin{align*}
a_0g(Y_0,Y_0)&=g(\mrm u, Y_0)=%
    g(\mrm u,-I_1\hat K_r-2vI_1K_0)=\tilde\om_1(K_r,\mrm u),\\%
a_1g(Y_1,Y_1)&=g(\mrm u, Y_1)=%
    g(\mrm u,\hat K_r+2vK_0)=2\xi(\mrm u)+\tilde g(K_r,\mrm u),\\%
a_2g(Y_2,Y_2)&=g(\mrm u, Y_2)=%
    g(\mrm u,I_3\hat K_r+2vI_3K_0)=-\tilde\om_3(K_r,u),\\%
a_3g(Y_3,Y_3)&=g(\mrm u, Y_3)=%
    g(\mrm u,-I_2\hat K_r-2vI_2K_0)=\tilde\om_2(K_r,u).%
\end{align*}
The expressions for $a_l$ become more compact in quaternionic notations. Indeed, if
we put $a=a_0+a_1i+a_2j+a_3k\in\mathbb H$ and recall the definition
(\ref{Eq_FormEta}) of 1--form $\eta$, then%
\begin{equation*}
a=\frac 1{4v+w}\bigl (2\xi(\mrm u)+\bar\eta(\mrm u)\bigr )i.
\end{equation*}

%Further, any tangent space of a \hK manifold carries an action of $\mathbb H$.
%In particular if $q\in\mathbb H$ we write $qY_0=\sum_{l=0}^3q_l Y_l$ for the
%sake of brevity.

Since $l_*=\duzhky{L_{\lambda(p)}}_* +d\lambda\, Y_0(l(p))$, we have $pr\, l_*\mrm
v=\duzhky{L_{\lambda(p)}}_*\mrm v'$ and therefore%
\begin{align*}
g\duzhky{pr\, l_*\mrm u\, , pr\, l_*\mrm v}&=%
    g\duzhky{\duzhky{L_{\lambda(p)}}_* \mrm u'\, , \duzhky{L_{\lambda(p)}}_*\mrm v'}\\%
&=\lambda^2 g\duzhky{\mrm u'\, , \mrm v'}\\%
&=\lambda^2 g\duzhky{\mrm u - aY_0 , \mrm v - bY_0}\\%
&=\lambda^2\bigl (g(\mrm u, \mrm v)- g(aY_0, \mrm v)-%
    g(\mrm u, bY_0)+g(aY_0, bY_0)\bigr ),
\end{align*}
where we also put $\mrm v=\mrm v'+bY_0$. Now it is easy to compute every single
summand in the last expression. Indeed, the first summand is given by formula
(\ref{MetricOnP}). Taking into account decompositions (\ref{Y_andK_r}) one obtains
$g(aY_0, \mrm v)= -\Re\duzhky{ai\duzhky{2\xi +\eta}(\mrm
v)}=(4v+w)^{-1}\Re\duzhky{\duzhky{2\xi +\bar\eta}(\mrm u)\duzhky{2\xi
+\eta}(\mrm v)}$. Since vectors $Y_l$ are pairwise orthogonal we get: % for the last term%
$g(aY_0,bY_0)=(4v+w)\Re\duzhky{a\bar b}=$\linebreak%
$(4v+w)^{-1}\Re\duzhky{\duzhky{2\xi+\bar\eta}(\mrm u)\duzhky{2\xi +\eta} (\mrm
v)}$. Finally, gathering all terms together one has after a simplification:%
\begin{equation*}
g_{\mt N}=\frac 1{4v+w} \duzhky{\tilde g +\frac 1v\xi^2-\frac
1{2(4v+w)}\duzhky{2\xi+\bar\eta}\odot\duzhky{2\xi+\eta}}.
%&=\frac 1{4v+w}\duzhky{\tilde g +\frac 1v\xi^2-\frac 12\psi^2}%
%    -\frac 1{2(4v+w)^2}\sum_{s=1}^3\duzhky{\imath_{K_r}\tilde\om_s}^2.
\end{equation*}

It is convenient to introduce a 1--form %
\begin{equation}\label{FormPsi}
\psi = \frac 1{g(Y_1,Y_1)}g(Y_1,\cdot)=%
\frac 1{4v+w}(2\xi + \imath_{K_r}\tilde g),
\end{equation}
which is a connection on the principal fibre bundle $P\2 N_0$ (assuming that $S_r^1$
acts freely on $P$),  i.e. it is $S_r^1$--invariant and $\psi (Y_1)=1$.
Then the expression for the metric takes the following form:%
\begin{equation}\label{MetricOnN}
g_{\mt N} =\frac 1{4v+w}\duzhky{\tilde g +\frac 1v\xi^2-\frac 12\psi^2}%
    -\frac 1{2(4v+w)^2}\sum_{l=1}^3\duzhky{\imath_{K_r}\tilde\om_l}^2.
\end{equation}

Arguments similar to those at the beginning of Section~\ref{Section_ReverseConstruction}
show that formula~(\ref{MetricOnN}) defines a metric on $N_0$.

%Arguing as at the beginning of Section~\ref{Section_ReverseConstruction}, one
%can show that formula~(\ref{MetricOnN}) defines a metric on $N_0$.

%In order to show that formula (\ref{MetricOnN}) defines a metric on $N_0$ one
%has to check that the bilinear form $g_{\mt N}$ is basic and invariant (\wrt to
%the action of $S_r^1$). A corresponding computation is similar to the one at
%the beginning of Section~\ref{Section_ReverseConstruction}.

%Checking that so defined symmetric tensor on $P$ is basic and invariant we
%obtain a metric on $N_0$.

The fundamental 4--form $\Om$ can be obtained in the similar manner. Indeed,%
\begin{equation*}
\begin{aligned}
\chi (\mrm u ,\mrm v)&=\om(pr\, l_*\mrm u , pr\, l_*\mrm v)%
    \lambda^2\om(\mrm u - aY_0,\mrm v - bY_0)=\\%
&=\lambda^2\duzhky{\om\duzhky{\mrm u,\mrm v} - \om\duzhky{aY_0,\mrm v}-\om (\mrm u, bY_0) +%
    \om(aY_0, bY_0)}.
\end{aligned}
\end{equation*}
Arguing similarly as we did when computing the metric, we obtain finally%
\begin{equation*}
\chi = \frac 1{4v+w}\tilde\om - %
    \frac 1{2(4v+w)^2}\duzhky{2\xi+\bar\eta}\wedge \duzhky{2\xi+\eta}.
\end{equation*}
Componentwise spelling of this formula is%
\begin{equation}\label{FormsChi}
\begin{gathered}
\chi_1=\frac 1{4v+w}\duzhky{\tilde\om_1-\psi\wedge\imath_{K_r}\tilde\om_1 }+%
    \frac 1{(4v+w)^2}\imath_{K_r}\tilde\om_2\wedge \imath_{K_r}\tilde\om_3,\\%
\chi_2=\frac 1{4v+w}\duzhky{\tilde\om_2-\psi\wedge\imath_{K_r}\tilde\om_2 } -%
    \frac 1{(4v+w)^2}\imath_{K_r}\tilde\om_1\wedge \imath_{K_r}\tilde\om_3,\\%
\chi_3=\frac 1{4v+w}\duzhky{\tilde\om_3-\psi\wedge\imath_{K_r}\tilde\om_3 }+%
    \frac 1{(4v+w)^2}\imath_{K_r}\tilde\om_1\wedge \imath_{K_r}\tilde\om_2.%
\end{gathered}
\end{equation}
With respect to the action of $S_r^1$ on $P$ all three forms
$\chi_l\in\Om^2(P)$ are basic, however only $\chi_1$ is invariant:%
\begin{align*}
\mathcal L_{Y_1}\chi_1 &= %
    \frac 1{(4v+w)^2}\duzhky{-2\imath_{K_r}\tilde\om_3\wedge \imath_{K_r}\tilde\om_3 +%
    \imath_{K_r}\tilde\om_2\wedge 2\imath_{K_r}\tilde\om_2}=0,\\%
\mathcal L_{Y_1}\chi_2 &=%
    \frac 1{4v+w}\duzhky{-2\tilde\om_3+\psi\wedge 2\imath_{K_r}\tilde\om_3 }+%
    \frac 1{(4v+w)^2} \imath_{K_r}\tilde\om_1\wedge 2\imath_{K_r}\tilde\om_2=%
    -2\chi_3,\\%
\mathcal L_{Y_1}\chi_3 &=%
    \frac 1{4v+w}\duzhky{2\tilde\om_2-\psi\wedge 2\imath_{K_r}\tilde\om_2 }-%
    \frac 1{(4v+w)^2} \imath_{K_r}\tilde\om_1\wedge 2\imath_{K_r}\tilde\om_3=%
    2\chi_2,%
\end{align*}
It follows that a 4--form%
\begin{equation}\label{FundamentalForm}
\Om=\chi_1\wedge\chi_1+ \chi_2\wedge\chi_2+ \chi_3\wedge\chi_3
\end{equation}
is basic and invariant and therefore descends to $N_0$. Integrability of such
defined \qK structure follows from integrability of the \hK structure on
$M_0=\mathcal H (\tilde M)=\mathcal U(N_0)$ \cite{Swann:91}.

\begin{thm}\label{Thm_QKmflds}
Let the assumptions of Theorem~\ref{Thm_HKmanifold} be satisfied. Then the metric
(\ref{MetricOnN}) and the fundamental 4--form (\ref{FundamentalForm}) define a \qK
structure on $N_0=\mathcal Q (\tilde M)=P/S_r^1$, where $\chi_l$ and $\psi$ are defined
by (\ref{FormsChi}) and (\ref{FormPsi}) respectively. Moreover $\mathcal Q (\tilde M)$
admits a \qK action of $S^1$ and its Swann bundle $\mathcal U(N_0)$ is $\mathcal H(\tilde
M)$.\qed
\end{thm}

%---END---\subsection{Quaternionic Flip.}

%---END---\section{S^1 Symmetry.}

\section{Examples}

\begin{ex}[$T^*\mathbb{CP}^n$ with the Calabi metric]\label{Example_CalabiMetric}
%\subsubsection{$T^*\mathbb{CP}^n$ with Calabi metric.}\label{Example_CalabiMetric}
The \hK quotient of $\mathbb H^{n+1}$ by $S^1$ acting by multiplication on the left with
respect to non--zero value of momentum map is topologically $T^*\mathbb{CP}^n$. Hitchin
\cite{Hitchin:86} showed, that the metric coincides with the one defined by Calabi
\cite{Calabi:79}. Therefore $\mathcal H(T^*\mathbb{CP}^n)=\mathbb H^{n+1}$ with its flat
metric and $\mathcal Q(T^*\mathbb{CP}^n)=\mathbb{HP}^n$ (in both cases with zero level
set of corresponding momentum map being removed).
\end{ex}

\begin{ex}[Flat manifold, adjoint action] Let us take a copy of quaternions
$\mathbb H_y$ as a manifold $\tilde M$ with the following action of $S_r^1:
(z,y)\mapsto zy\bar z$ (one can also regard $\mathbb H$ as $T^*\mathbb C$ with
fibrewise action of $S_r^1$; see also remark \ref{Rem_CotangentBundle}). In this
case $1/4(w+2\tilde\rho) = 1/2(y_2^2+y_3^2)$, where
$y=y_0+y_1i+y_2j+y_3k$. Adding $1/2$ we may write function $v$ in the form%
\begin{equation*}
v=\frac 12(1-y_2^2-y_3^2)
\end{equation*}
and it is positive on $\mathbb R_{y_0\, y_1}^2\times \mrm D_{y_2\, y_3}^2,$ where
$\mrm D^2\subset\mathbb R^2$ is an open disc of radius~$1$. The principal bundle $P$
is trivial and therefore $\mathcal Q(\mathbb R^2\times \mrm D^2)=\mathbb R^2\times
\mrm D^2$
with the following metric:%
\begin{align*}
g_{\mt N}=\frac 1{2(1+d)}&\left ( %
\frac {1-d}{1+d}\,\Re dy\otimes d\bar y  +%
\frac{4d}{1-d^{\, 2}}\duzhky{y_0dy_1 +y_3dy_2}^2 \right .\\ %
&\qquad -\left .\frac 1{1+d}\duzhky{y_0dy_1 +y_3dy_2}\odot \duzhky{y_2dy_3 -
y_3dy_2} \right ),
\end{align*}
where $d=y_2^2 +y_3^2$. Therefore the above metric is Einstein and self--dual.
However it is incomplete.

Similarly, one can compute the metric and symplectic forms on $\mathcal H (\mathbb
R^2\times \mrm D^2)=\Hstar\times\mathbb R^2\times \mrm D^2$ but the metric is also
incomplete.
\end{ex}
%===============\begin{ex}\textbf{Flat manifold, adjoint action.}

\begin{rem}\label{Rem_CotangentBundle}
The above manifolds are examples of a large class of \hK manifolds admitting
$S_r^1$--action. Namely Kaledin \cite{Kaledin:97} and independently Feix \cite{Feix:01}
proved the existence of a \hK metric on a (neighborhood of the zero section of) cotangent
bundle $T^*\mathcal M$ to a real--analytic K\"ahler manifold $\mathcal M$. The above
examples show that the function $v$ can be both positive everywhere and only on a proper
open subset of $T^*\mathcal M$.
\end{rem}

\begin{ex}[Gibbons--Hawking Spaces]\label{Ex_GibbonsHawking}
All \hK four-manifolds with $S^1$--symmetry were described by Gibbons and Hawking
\cite{GibbonsHawking:78} and their construction is as follows. If $Z^4$ is \hK and admits
$S^1$--symmetry with a \Kvf $K_0$, then its \hK momentum map $\mu=\mu_1i+\mu_2j+\mu_3k$
represents $Z$ as a fibration over $\mathbb R^3$ with generic fibre $S^1$, so that,
excluding critical points of the momentum map, one can write the metric as%
\begin{equation}\label{Eq_GibbonsHawkingMetric}
g_{\mt{GH}}=\nu\,\duzhky{dx_1^2+ dx_2^2+dx_3^2}+\nu^{-1}\xi^2,\qquad x_l=\mu_l,\
l=1,2,3,
\end{equation}
where $\nu:\mathbb R^3\2 \mathbb R_{>0},\ \nu^{-1} =\| K_0 \|^{2}$, and $\xi$ is a
connection form. It is then an easy exercise to write down 2--forms which are closed
provided
\begin{equation}\label{Eq_GH}
F_\xi = -* d\nu.
\end{equation}
It follows form the Bianchi identity that $\nu$ is harmonic. We would like to point
out that $Z^4$ is determined by the function $\nu$ (harmonic and positive) since the
connection $\xi$ can be found from equation (\ref{Eq_GH}).
The Gibbons--Hawking--Ansatz is a choice of a particular function $\nu:$%
\begin{equation*}
\nu(x)=\sum\limits_{i=1}^n \frac 1{|\, x-y_i|},\qquad y_i\in \mathbb R^3.
\end{equation*}
In general the above four-manifold does not admit an $S_r^1$--action. However when all
the poles $y_i$ of the function $\nu$ lie on one line (say $x_1$--axis), then such action
does exist; its projection to $\mathbb R^3\cong \Im\mathbb H$ is then $(z,x)\mapsto
zx\bar z,\ x\in \Im\mathbb H$. A direct (and tedious) computation shows that the function
$v$ is positive everywhere and therefore the construction $\mathcal H(Z)$ is defined on
the whole Gibbons--Hawking space $Z$. Alternatively, one can observe \cite{Galicki:87}
that Gibbons--Hawking spaces can be obtained as \hK reductions of a flat space acted upon
by a torus \wrt a nonzero value of the corresponding momentum map, i.e. Gibbons--Hawking
spaces are examples of toric \hK manifolds~\cite{BielawskiDancer:00}. Therefore, if the
value of the momentum map is chosen properly, the action of $S^1_r$ can be obtained from
the corresponding action on the flat space. In this case one can also show that the
manifold $\mathcal H(Z)$ is also toric.

\end{ex}
%================\begin{ex}[Gibbons--Hawking Spaces]

%---END---\section{Examples.}

\section{Indeterminacy of function $v$: further examples}\label{Section_IndeterminacyOfV}

As we have seen in Section \ref{Section_ReverseConstruction} the function $v$ is defined
by formula (\ref{Function_v}) only up to a constant and this has strong consequences as
we will see below. Recall that the action of $S_r^1$ should be lifted from $\tilde M$ to
$P$ such that its \Kvf $Y_1$ equals to $\hat K_r+ 2vK_0$ (see (\ref{Y_andK_r})). This
implies that we are free to take $\tilde v=v+m/2$ instead of $v$, where $m$ is an
integer, whenever $v+m/2$ remains everywhere positive. However in this case one needs to
modify the lifting of the $S_r^1$--action to $P$ to get that the \Kvf is given by $\hat
K_r+ 2\tilde vK_0$. Therefore we get that the manifolds $\Hstar\times_{S_r^1}P$ and
$P/S_r^1$ again carry \hK and \qK structures correspondingly, where the modified action
of $S_r^1$ is implied. It turns out that the modification of the $S_r^1$--action can
change the topology of the $\mathcal H$ and $\mathcal Q$ constructions.

In the rest of this section we carry out the above observation in details for the case of
cotangent bundle of a complex Grassmannian. Hyper\kahler metrics on cotangent bundles of
Grassmannians with required symmetries were obtained long ago (see \cite{Galicki:87} and
references therein). More generally, Nakajima~\cite{Nakajima:94} constructed such metrics
on cotangent bundles of quiver varieties, however we shall consider only Grassmannians
for the sake of simplicity. First we review construction of \hK structure on
$T^*Gr_k(\mathbb C^n)$ and then illustrate the impact of the modification of
$S^1_r$--action.

Choose the left quaternionic structure on the flat space $M_{n,k}(\mathbb H)$ consisting
of matrices with $n$ rows and $k$ columns. We have
\begin{equation*}
g(A_1, A_2)=\mathrm{Re\, tr} (A_1 \bar A_2^t),\quad %
\om(A_1, A_2)= \mathrm{Im\, tr} (A_1 \bar A_2^t).
\end{equation*}
The group $U(k)$ acts on $M_{n,k}(\mathbb H)$ by multiplication on the right. Write
$A=B+C^tj$, where $B\in M_{n,k}(\mathbb C)$ and $C\in M_{k,n}(\mathbb C)$. Then the \hK
momentum map $\mu=\mu_{\ssst{\mathbb R}}i+\mu_{\ssst{\mathbb C}}j$ is given by%
\begin{equation*}
\mu_{\ssst{\mathbb R}}(B,C)=\frac 12 \left (\bar B^tB-C\bar C^t \right),\quad %
\mu_{\ssst{\mathbb C}}(B, C)=-CB.
\end{equation*}

Probably the easiest way to see that the \hK reduction of $M_{n,k}(\mathbb H)$ is
isomorphic to $T^*Gr_k(\mathbb C^n)$ is to observe~\cite{Nakajima:94} that%
\begin{equation}\label{Eq_TGrAsQuotient}
\begin{aligned}
\mu^{-1}(i)&/U(k)\cong\\%
&\{ (B,C)\in M_{n,k}(\mathbb C)\times M_{k,n}(\mathbb C)\, |\ \mathrm{rk}\, B=k, BC=0
\}/GL_k(\mathbb C),
\end{aligned}
\end{equation}
where $GL_k(\mathbb C)$ acts on pairs of matrices as follows: %
$(B, C)\cdot g=(Bg, g^{-1}C)$. Think about $B$ as a k-frame in $\mathbb C^n$. Then %
$\{ B\, |\ \mathrm{rk}\, B=k\}/GL_k(\mathbb C)= Gr_k(\mathbb C^n)$ and it follows that
the right-hand side of (\ref{Eq_TGrAsQuotient}) is $S\otimes Q^\vee$, where $S$ and $Q$
denote tautological and quotient vector bundles over $Gr_k(\mathbb C^n)$ correspondingly.
Recalling that $TGr_k(\mathbb C^n)\cong S^\vee\otimes Q$ we get the result.

The action of $S^1_r$, inerited from the permuting action of $\Hstar$ on $M_{n,k}(\mathbb
H)$, is the following one\footnote{notice that this action is different from the one
considered by Nakajima in~\cite{Nakajima:94}}: $z\cdot [B, C]=[zB, zC]$. For such action
the function $v$ defined by (\ref{Function_v}) must be positive everywhere on
$T^*Gr_k(\mathbb C^n)$. Indeed, this follows from the following observation. The space
$M_{n,k}/\hskip-0.1cm/\hskip-0.1cm/_{\mu=i}U(k)$ can be obtained in two steps: first
consider $M_{n,k}/\hskip-0.1cm/\hskip-0.1cm/_{\mu=0}SU(k)$ and then take its \hK
reduction \wrt $S^1\subset U(k)$. Since the space
$M_{n,k}/\hskip-0.1cm/\hskip-0.1cm/_{\mu=0}SU(k)$ inherits permuting action of $\Hstar$,
the statement follows. Notice also that the corresponding $S^1$--principal bundle
$P=\mu^{-1}_{S^1}(i)\subset M_{n,k}/\hskip-0.1cm/\hskip-0.1cm/_{\mu=0}SU(k)$ is pull-back
of the principal bundle of $\Lambda^{top}S$.

Take $\tilde v=v+m/2$ instead of $v$, where $m$ is a positive integer. Then the modified
action of $S^1_r$ on $\mu^{-1}(i)$ is given by%
\begin{equation}\label{Eq_ModifiedActionOfS1}
z\cdot (B, C)=(Bz^{m+1}, z^{1-m}C).
\end{equation}
Before proceeding we need the following Lemma.

\begin{lem}\label{Lemma_S1Bundles}
Let $P\ra X$ be an $S^1$--principal bundle and $L\ra X$ be the corresponding line bundle.
Consider the following action of $S^1$ on $P\times\mathbb C:\ z\cdot
(p, w)=(pz^r, z^sw),$ where $r$ and $s$ are integers and $r$ is positive. Then%
\begin{equation*}
(P\times\mathbb C )/S^1\cong L^{-s}.
\end{equation*}
\end{lem}

\begin{proof}
Let $Q_{r,\, s}$ denote the space $P\times\mathbb C$ with the action of $S^1$ as in the
statement of the lemma. Then we have an equivariant map $Q_{r,\, s}\ra Q_{r,\, rs},\
(p,w)\mapsto (p,w^r)$. Clearly it is surjective; although it is not injective, it
descends to a bijective map of quotients $Q_{r,\, s}/S^1\ra Q_{r,\, rs}/S^1$. But the
last quotient is exactly $L^{-s}$.\qed
\end{proof}

Therefore the action (\ref{Eq_ModifiedActionOfS1}) can be replaced by the following one%
\begin{equation*}
z\cdot (B, C)=(Bz, z^{-1}Cz^{2-m}).
\end{equation*}
This action is induced by the inclusion $S^1\subset U(k)$ followed by the action of
$U(k)$%
\begin{equation*}
(B,C)\cdot g=(Bg, g^{-1} C(\det g)^{-s}),
\end{equation*}
provided $m=ks+2,\ s\in\mathbb Z$. If $L$ denotes the top exterior power of the tautological bundle of %
$Gr_k(\mathbb C^n)$, then we get%
\begin{equation*}
\mu^{-1}_{S^1}(i)/S^1_r\cong \mu^{-1}(i)/U(k)\cong L^s\otimes T^*Gr_k(\mathbb C^n).
\end{equation*}
Thus, summing up we get the following result.

\begin{thm}
The total space of $L^s\otimes T^*Gr_k(\mathbb C^n),\ s\in\mathbb Z,\ ks+2\ge 0$ carries
a \qK structure with positive scalar curvature. Its Swann bundle is the total space of
$L\oplus L^{-1}\oplus L^s\otimes T^*Gr_k(\mathbb C^n)$.\qed
\end{thm}

In case of $k=1, n=2$ one obtains \qK structures on total spaces of $\mathcal O_{\mathbb
P^1}(-s-2)=\mathcal O_{\mathbb P^1}(-m)$ for $m\ge 0$, which is certainly diffeomorphic
to $\mathcal O_{\mathbb P^1}(m)$. These spaces were obtained by Galicki and
Lawson~\cite{GalickiLawson:88} as open subsets of four-dimensional \qK orbifolds.

%---END---\section{Indeterminacy of function $v$: further examples}\label{Section_IndeterminacyOfV}

\section{K{\" a}hler structure on $N_0$}

In contrast to a \hK manifold, almost complex structures of a \qK manifold $N$ are
defined only locally, i.e. we have a distinguished rank 3 subbbundle
$\bundlei\subset End(TN)$ called a \textit{structure bundle}\index{Structure
bundle}, which locally admits a basis consisting of three almost complex structures
with quaternionic relations. Since the metric induces an isomorphism $TN\cong T^*N$,
one gets an embedding of $\bundlei$ in $\Lambda^2T^*N$. Locally this is given by
passing from an almost complex structure $I$ to the associated 2--form $\om_I(\cdot
,\cdot)=g(\cdot , I\cdot)$. We will not distinguish between $\bundlei$ and its image
in $\Lambda^2T^*N$. An analogue of a momentum map can be defined in the \qK context,
but now it will be a section of a structure bundle (see \cite{Galicki:87} for
details).

\begin{thm}\label{Thm_KahlerStructure}
Let $N$ be a \qK manifold of positive scalar curvature. Suppose also that $N$ admits
a \qK action of $S^1$ with momentum section $\mu_{\ssst N} \in\Gamma(\bundlei )$.
Then $N_0=N\setminus\{\mu_{\ssst N} =0 \}$ is a K{\" a}hler manifold.
\end{thm}
\begin{proof}
Let $M$ be the Swann bundle of $N$. Then $M$ admits a \hK action of $S^1$
\cite{Swann:91}. Let $\mu = \mu_1i +\mu_c j: M\ra \Im\mathbb H$ be its momentum map.
Since the function $\mu_c: M\ra \mathbb C$ is $I_1$--holomorphic, the submanifold $M_c=\{
m\in M : \mu_c(m)=0\}$ has an induced \kahler structure. Further, $M_c^+=\{ m\in M_c:
\mu_1(m)>0\}$ is an open submanifold of $M_c$. The group $S_r^1$ preserves $I_1$ and one
may consider the \kahler reduction of $M_c^+$ \wrt a non--zero value of the momentum map
$:\ M_c^+/\hskip-0.1cm/S_r^1\cong M_c^+/\mathbb C_r^*$. It remains to observe that
$M_c^+=\mu^{-1}(i)\times\mathbb R_{>0}=P\times\mathbb R_{>0}$ and therefore
$M_c^+/\mathbb C_r^*\cong P/S_r^1\cong N_0$.\qed
\end{proof}
When a \qK manifold $N$ admits an action of $S^1$, one can normalize the momentum
section $\mu_{\ssst N} \in\Gamma(\bundlei )$ and consider it as an almost complex
structure $\hat I$ over $N_0$. It turns out that $\hat I$ is
\textit{integrable}~\cite{Battaglia:99, Salamon:99} and it is easy to see from the
proof that it coincides with the complex structure implied by
Theorem~\ref{Thm_KahlerStructure} (our proof of the above theorem itself represents
an alternative proof of the integrability of $\hat I$ in case when $N$ has positive
scalar curvature). Although the complex structure $\hat I$ is a section of the
structure bundle $\bundlei$, the \textit{\kahler} metric of $N_0$ must not coincide
with the \qK one as we will see in the sequel. Note also that $\bundlei$ does not
admit a section which defines an integrable complex structure on the whole manifold
$N$ (see \cite{Alekseevsky:98} for extensive discussion of this phenomenon). Taking
this into account, one may consider $N_0$ as "the largest" open submanifold of $N$
where it is still possible to choose an integrable complex structure.
%\medskip

Our next aim is to express the \kahler structure of $N_0$ similarly to the \qK one
(see section \ref{Section_QuaternionicFlip}).

Recall that $N_0\cong P/S_r^1=M_c^+/\mathbb C_r^*$. In order to get a metric and
\kahler form on $N_0$ we have to express $N_0$ as a \kahler reduction, i.e. we have
to fix a level set of momentum map and divide it by $S_r^1\subset \mathbb C_r^*$. In
our case the momentum map of the $S_r^1$--action is nothing else but the \hK
potential $\rho$ (restricted to $M_c^+$). Recall also that we have an isomorphism
(\ref{Eq_Isomorphism_l}) between $P$ and $Q=\rho^{-1}(-1/2)\cap M_c^+$. Further, one
has $\Tdot M_c= \lspan (I_2K_0, I_3K_0)^\perp\subset \Tdot M$ and $\Tdot
\rho^{-1}(-1/2)=\lspan (Y_0)^\perp$. It follows that $\Tdot Q=\lspan (I_2K_0,
I_3K_0, Y_0)^\perp$ because $Y_0$ is perpendicular to both $I_2K_0$ and $I_3K_0$
(see (\ref{Y_andK_r})). In particular $Y_1\in TQ$; this also follows from the fact
that $S_r^1$ preserves $Q$. Further, the \kahler reduction procedure implies that
$\Tdot N_0$ is identified with $\lspan (Y_1)^\perp\subset \Tdot Q$ and the \kahler
form and metric are obtained as a restriction of the corresponding tensors to
$\lspan (Y_1)^\perp$. Remark that the \qK metric was obtained as the one induced on
a different subbundle, namely on $\lspan (Y_1, Y_2, Y_3)^\perp\subset TQ$.

Let $\mrm u\in T_pP$. Then we may decompose $\mrm u = \mrm u' +\psi (\mrm u)Y_1$,
where $\mrm u'$ is orthogonal to $Y_1$. Now denote by $\Pi$ an orthogonal projector
on $\lspan (Y_1)^\perp$ in $TQ$. Then for the \kahler
metric $\hat g_{\mt N}$ we have:%
\begin{align*}
\hat g_{\mt N}&(\mrm u , \mrm v)= g\bigl (\Pi\, l_*\mrm u, \Pi\, l_*\mrm v\bigr )\\%
    &= g\duzhky{\duzhky{L_{\lambda(p)}}_* \mrm u' + d\lambda (\mrm u)Y_0(\lambda(p)p),\ %
            \duzhky{L_{\lambda(p)}}_* \mrm v' + d\lambda (\mrm v)Y_0(\lambda(p)p)}\\%
    &= \lambda^2g\bigl (\mrm u - \psi(\mrm u)Y_1 +d\lambda(\mrm u)Y_0\, ,\ \mrm v - %
            \psi(\mrm v)Y_1 +d\lambda(\mrm v)Y_0\bigr )%
\end{align*}
\begin{align*}
    &= \lambda^2\Bigl ( g\duzhky{\mrm u ,\mrm v}- \psi(\mrm u)g\duzhky{Y_1,\mrm v}-%
            \psi(\mrm v)g\duzhky{Y_1,\mrm u}+d\lambda\duzhky{\mrm u}g\duzhky{Y_0,\mrm v} %
            d\lambda\duzhky{\mrm v}g\duzhky{Y_0,\mrm u} \Bigr ).
\end{align*}

As we already know $g\duzhky{\mrm u ,\mrm v}=\bigl ( \tilde g +v^{-1}\xi^2 \bigr ) (\mrm
u,\mrm v)$. By the definition of $\psi$ one has $g(Y_1, \mrm v)=
(4v+w)\psi(\mrm v)$. Further, %
$g(Y_0, \mrm v)= g\bigl (-I_1K_r-2vI_1K_0,\,\hat{\mrm v}
+\xi(\mrm v)K_0\bigr )= \imath_{K_r}\tilde\om_1 (\mrm v)$. Therefore we obtain %
$\hat g_{\mt N}=\lambda^2\bigl ( \tilde g+ v^{-1}\xi^2 - (4v+w)\psi^2 +
d\lambda\odot\imath_{K_r}\tilde\om_1 \bigr )$. Since $d\lambda =
(4v+w)^{-3/2}\imath_{K_r}\tilde\om_1$ we may finally write%
\begin{equation*}
\hat g_{\mt N}=\frac 1{4v+w}\, \tilde g +\frac 1{v(4v+w)}\, \xi^2 -\psi^2 + \frac
1{(4v+w)^{5/2}}\, \duzhky{\imath_{K_r}\tilde\om_1}^2 .
\end{equation*}

The \kahler form $\hat\om_{\ssst N}$ may be obtained in a similar manner. Indeed, %
\begin{align*}
\hat\om_{\ssst N} (\mrm u, \mrm v)&= \om_1 (\Pi\, l_*\mrm u, \Pi\, l_*\mrm v) \\%
    &=\om_1\Bigl (\duzhky{L_{\lambda(p)}}_* \mrm u' + d\lambda (\mrm u)Y_0\bigl (l(p)\bigr ),\ %
            \duzhky{L_{\lambda(p)}}_* \mrm v' + d\lambda (\mrm v)Y_0\bigl (l(p)\bigr )\Bigr )\\%
    &=\om_1\duzhky{\duzhky{L_{\lambda(p)}}_* \mrm u',\ %
            \duzhky{L_{\lambda(p)}}_* \mrm v' }\\%
    &=\lambda^2\om_1\bigl (\mrm u - \psi(\mrm u)Y_1,\ \mrm v - \psi(\mrm v)Y_1
    \bigr).
\end{align*}
Since $\om_1(\mrm u, \mrm v)=\tilde\om_1(\mrm u, \mrm v)$ and $\om_1(Y_1, \mrm
u)=g(\hat K_r +2vK_0, I_1\hat{\mrm u} +\xi(\mrm
u)I_1K_0)=\imath_{K_r}\tilde\om_1 (\mrm u)$, we obtain the \kahler form as%
\begin{equation*}
\hat\om_{\ssst N} = \frac 1{4v+w}\, \Bigl ( \tilde\om_1 -
\psi\wedge\imath_{K_r}\tilde\om_1 \Bigr ).
\end{equation*}

\begin{rem}
As we have already remarked, we may regard the form $\psi$ as a connection on
the $S_r^1$--principal bundle $P\ra N_0$. Lets compute its curvature. We have%
\begin{equation*}
F_\psi= -\frac 1{(4v+w)^2}(4dv+dw)\wedge (2\xi +\imath_{K_r}\tilde g) + %
    \frac 1{4v+w} \duzhky{2d\xi + d\, \imath_{K_r}\tilde g\, }.
\end{equation*}
It follows from equations (\ref{Eq_F}) and  (\ref{Eq_dv}) that%
\begin{equation*}
F_\psi=-\frac 2{4v+w}\bigl ( \imath_{K_r}\tilde\om_1\wedge\psi +\tilde\om_1 \bigr )
= -2\hat\om_{\ssst N}
\end{equation*}

This observation provides an "intrinsic" interpretation of the \kahler form
$\hat\om_{\ssst N}$ in the following sense. Let $N$ be a \qK manifold with positive
scalar curvature and $F\ra N$ be the principal $SO(3)$ bundle associated to the
structure bundle $\bundlei$. Observe that $F$ is equipped with the natural
connection induced by the Levi--Civita one. Suppose also that $N$ admits a \qK
action of the circle and denote by $\mu_{\ssst N}$ its momentum section. As it was
explained above one can think about $\mu_{\ssst N}$ on $N_0=N\setminus\{\mu_{\ssst
N} =0 \}$ as a section of $F$ (restricted to $N_0$). This means that we get
$S^1$--subbundle $Q$ of $F$. The curvature of the induced connection $\psi$ is a
$-1/2$--multiple of the \kahler form $\hat\om_{\ssst N}$.
\end{rem}

\medskip
\textbf{Acknowledgement.} This paper is based upon part of the author's doctoral
thesis. I would like to thank my supervisor Prof.\,V.\,Pid\-stry\-gach for his
constant encouragement and support.

\end{document}